\documentclass[10pt]{article}

\usepackage{amsfonts,amsmath,graphicx}
\usepackage[pdfpagemode=None,colorlinks=true,linkcolor=blue,citecolor=blue]{hyperref}

\newtheorem{proposition}{Proposition}
\newtheorem{remark}{Remark}
\newtheorem{theorem}{Theorem}

\newcommand{\assign}{:=}
\newenvironment{itemizedot}
  {\begin{itemize}
  }{\end{itemize}}

\newcommand{\Z}[1][]{\ensuremath{{\mathbb{Z}^{#1}} }}
\newcommand{\C}[1][]{\ensuremath{{\mathbb{C}^{#1}} }}
\newcommand{\R}[1][]{\ensuremath{{\mathbb{R}^{#1}} }}
\renewcommand{\H}[1][]{\ensuremath{{\mathbb{H}^{#1}} }}
\renewcommand{\S}[1][]{\ensuremath{{\mathbb{S}^{#1}} }}

\def\Re{ \mathrm{Re}\, }
\def\Im{ \mathrm{Im}\, }
\newcommand{\<}{\langle}
\renewcommand{\>}{\rangle}

\date{}
\title{Cyclic and ruled Lagrangian surfaces in complex Euclidean space}
\author{Henri Anciaux and Pascal Romon}

\begin{document}
\maketitle

\medskip

\begin{abstract}
We study those Lagrangian surfaces in complex Euclidean space
which are foliated by circles or by straight lines. The former, 
which we call \emph{cyclic}, come in three types, each
one being described by means of, respectively, a planar curve, a
Legendrian curve in the 3-sphere or a Legendrian curve in the anti-de Sitter 
3-space. We describe ruled Lagrangian surfaces
and characterize the cyclic and ruled Lagrangian surfaces
which are solutions to the self-similar equation of the Mean
Curvature Flow. Finally, we give a partial result in the case of
Hamiltonian stationary cyclic surfaces.

\end{abstract}

\section*{Introduction}

In this paper, we classify the Lagrangian surfaces of $\mathbb{C}^2$ which
are foliated either by round circles (henceforth called \emph{cyclic}
surfaces) or by straight lines (\emph{ruled} surfaces). This completes a
former paper of the authors together with Ildefonso Castro \cite{ACR} in which all
Lagrangian submanifolds of $\mathbb{R}^{2 n} \simeq \mathbb{C}^n$, with $n
\geq 3$, which are foliated by round $( n - 1 )$-spheres were characterized.
The reason for the lower bound on the dimension was the following: since the
submanifold is Lagrangian, any spherical leaf must be isotropic; when the
dimension of this leaf is at least two, it spans a linear space which is
itself Lagrangian. This observation simplifies the structure of the problem,
roughly speaking by reducing the underlying group structure from $SO ( 2 n )$
to $U ( n )$. However this reduction no longer holds in dimension two, see for
instance the Lagrangian cylinder $ \S^1 \times \mathcal{L}$ where $\mathcal{L}$ 
is a real line of $\C$: this
(Lagrangian) surface is foliated by circles which are contained in
non-Lagrangian (actually complex) planes. Other examples are the Hopf tori
studied by Pinkall in \cite{P}. As expected the situation is richer in dimension
two, and actually cyclic Lagrangian surfaces come in three families, each
one being described by means of, respectively, a planar curve, a
Legendrian curve in the 3-sphere or a Legendrian curve in the anti
de Sitter 3-space  (Theorems 1 and 2). In the following, they will be 
denoted as  type I, II and III surfaces.

In Section 1 we classify the cyclic Lagrangian
surfaces when all the centers of the circles coincide. We call those surfaces
\emph{centered} cyclic. In Section 2 we treat the general case which amounts
to adding a convenient translation term. In Section 3 we apply this
characterization to finding self-similar cyclic surfaces, that is
those surfaces which are solutions of the following elliptic PDE:
$$H + \lambda X^{\perp} = 0,$$
where $H$ denotes the mean curvature vector of the surface and $X^\perp$ the
normal component of its position vector.
 The case of positive (resp. negative) $\lambda$ corresponds
to the case of a self-shrinking (resp. self-expanding) soliton of the Mean Curvature
Flow (see \cite{A}).
 We show that a self-similar Lagrangian cyclic surface
is either a centered surface of type I
 as described in \cite{A} or the Cartesian product $ \S^1(r) \times \Gamma$ of a circle $\S^1(r)$
 with a planar self-shrinking curve $\Gamma$. Such curves have been studied in detail
 in \cite{AL}.

 Section 4 is devoted to the Hamiltonian stationary equation.
A Lagrangian surface is said to be \em Hamiltonian stationary \em
if its area is critical for compactly supported Hamiltonian
variations. Such a surface is characterized by the fact that its
Lagrangian angle $\beta$ is harmonic with respect to the induced
metric (cf Section 4 for more details). We prove that both in
types I and II cases, Hamiltonian stationary surfaces must be
centered, and we describe them. The study of the type III appears
to be extremely difficult to handle by manual computation, however
we conjecture that Hamiltonian stationary type III surfaces are again centered. 
Examples of such surfaces are described explicitly in \cite{CC}.

Finally we give in the last Section a description of ruled
Lagrangian surfaces (notably self-similar ones) using an analogous method,
recovering more simply a known result from Blair \cite{B}.

\medskip

Acknowledgment: the Authors wish to thank Ildefonso Castro for pointing out a
better approach to type III surfaces due to Chen and himself \cite{CC}.

\bigskip

\section{Centered cyclic Lagrangian surfaces}

Let $\Sigma$ a surface of $\mathbb{R}^4$ foliated by circles with common
center located at the origin of $\mathbb{R}^4$. Locally, $\Sigma$ may be
parametrized by the following immersion:
\[ \begin{array}{cccl}
     X : & I \times \mathbb{R} / 2 \pi \mathbb{Z} & \to & \mathbb{R}^4\\
     & ( s, t ) & \mapsto & r ( s ) ( e_1 ( s ) \cos t + e_2 ( s ) \sin t ),
   \end{array} \]
where $r ( s )$ is a positive function and $( e_1 ( s ), e_2 ( s ) )$ is an
orthonormal basis of the plane containing our circle.

From now on we shall assume that $\Sigma$ is Lagrangian with respect with some
complex structure $J$. We will often identify $\mathbb{R}^4$ with
$\mathbb{C}^2$ in such a way that $J$ is the complex multiplication by $i$.
Denote by $K \assign \langle e_1, J e_2 \rangle = - \langle e_1, J e_2 \rangle$
the \textit{K\"{a}hler angle} of the plane $e_1 \wedge e_2 $. Note that the
vanishing of $K$ means that $e_1 \wedge e_2$ is also Lagrangian. In this case
the analysis done in \cite{ACR} holds and $X$ takes the following form: $X ( s, t )
= r ( s ) e^{i \phi ( s )} ( \cos t, \sin t )$ making use of the above
identification.

Denoting by subscripts the partial derivatives (the prime corresponding also
to the derivative in $s$ for functions of just one variable), the Lagrangian
assumption is equivalent to $\langle X_s, J X_t \rangle = 0$. Since
\[ X_s = r' ( e_1 \cos t + e_2 \sin t ) + r ( e'_1 \cos t + e'_2 \sin t ) \]
\[ X_t = r ( e_2 \cos t - e_1 \sin t ), \]
we see that
\begin{eqnarray*}
  \langle X_s, J X_t \rangle & = & r r' ( K \cos^2 t + K \sin^2 t )
+ r^2  \left( \cos^2 t \langle e'_1, J e_2 \rangle - \sin^2 t \langle e'_2, J e_1 \rangle  \right)\\
 & & + r^2 \cos t \sin t \left( \langle e'_2, J e_2 \rangle -
  \langle e'_1, J e_1 \rangle \right) \\
  & = & r r' K + \frac{r^2}{2}  \left( \left\langle e_1', J e_2 \right\rangle -
  \langle e'_2, J e_1 \rangle \right) \\
&&+ \frac{r^2 \cos 2 t}{2}  \left( \langle e'_1, J e_2 \rangle + \langle  e'_2, J e_1 \rangle \right) \\
  & & + \frac{r^2 \sin 2 t}{2}  \left( \langle e'_2,
  J e_2 \rangle - \langle e'_1, J e_1 \rangle \right)\\
  & = & \frac{2 r r' K + r^2 K'}{2} + \frac{r^2 \cos 2 t}{2}  \left( \langle
  e'_1, J e_2 \rangle + \langle e'_2, J e_1 \rangle \right) \\
  & & + \frac{r^2 \sin 2
  t}{2}  \left( \langle e'_2, J e_2 \rangle - \langle e'_1, J e_1 \rangle
  \right)
\end{eqnarray*}
which holds for all $t$.
The vanishing of the constant term implies that $r^2 K = C$ for some real
constant $C$. If this constant vanishes, we recover the case $K = 0$ mentioned
above and treated in \cite{ACR} (cf also \cite{A}), so we may assume that $C \neq 0$.
Thus $r$ is completely determined by $K$, and both are non zero.

The two remaining conditions are
\begin{equation}
  \langle e'_1, J e_2 \rangle = - \langle e'_2, J e_1 \rangle \;, \quad \langle
  e'_1, J e_1 \rangle = \langle e'_2, J e_2 \rangle .
  \label{remainingconditions}
\end{equation}
In order to make sense of these, we will identify $\mathbb{R}^4$ with
$\mathbb{H}$, in such a way that the complex structure is given by the left
multiplication by the quaternion $i$. Then any element in $S O ( 4 )$ can be
written as $x \mapsto p x q^{- 1}$ where $p, q$ are two unit quaternions.
Notice that right multiplication by $q^{-1}$ corresponds exactly to the elements of
$S U(2)$.
Since $S O ( 4 )$ acts transitively on pairs of orthonormal vectors, we may write
$e_1$ and $e_2$ as the respective images of $1$ and $i$, so that $e_1 = p q^{-
1}$ and $e_2 = p i q^{- 1}$. Note that $( p, q )$ is not uniquely determined;
rather we have a gauge freedom by right multiplication by $e^{i \theta}$ on $(
p, q )$. Finally, we may assume if needed that $q(0)$ takes any prescribed value,
since we consider surfaces up to $U(2)$ congruence.

Then the conditions in (\ref{remainingconditions}) read as
\[ \left\{\begin{array}{l}
     \left\langle p' q^{- 1} - p q^{- 1} q' q^{- 1}, i p i q^{- 1} \right\rangle +
     \left\langle p' i q^{- 1} - p i q^{- 1} q' q^{- 1}, i p q^{- 1} \right\rangle
     = 0\\
     \left\langle p' q^{- 1} - p q^{- 1} q' q^{- 1}, i p q p q^{- 1} \right\rangle -
     \left\langle p' i q^{- 1} - p i q^{- 1} q' q^{- 1}, i p i q^{- 1} \right\rangle
     = 0
   \end{array}\right. \]
so, multiplying on the left by $p^{- 1}$ and the right by $q$, and further by
$i$ on the left in the second bracket,
\[ \left\{\begin{array}{l}
     \left\langle p^{- 1} p' - q^{- 1} q', p^{- 1} i p i \right\rangle -
     \left\langle p^{- 1} p' + i q^{- 1} q' i, p^{- 1} i p i \right\rangle = 0\\
     \left\langle p^{- 1} p' - q^{- 1} q', p^{- 1} i p \right\rangle -
     \left\langle p^{- 1} p' + i q^{- 1} q' i, p^{- 1} i p \right\rangle = 0
   \end{array}\right. \]
that is
\[ \left\{\begin{array}{l}
     \left\langle q^{- 1} q' + i q^{- 1} q' i, p^{- 1} i p i \right\rangle = 0\\
     \left\langle q^{- 1} q' + i q^{- 1} q' i, p^{- 1} i p \right\rangle = 0 .
   \end{array}\right. \]
Splitting and multiplying left and right by $i$ in the second bracket yields
\[ \left\{\begin{array}{l}
     \left\langle q^{- 1} q', p^{- 1} i p i \right\rangle - \left\langle q^{- 1}
     q', i p^{- 1} i p \right\rangle = 0\\
     \left\langle q^{- 1} q', p^{- 1} i p \right\rangle + \left\langle q^{- 1}
     q', i p^{- 1} i p i \right\rangle = 0
   \end{array}\right. \]
\[ \left\{\begin{array}{l}
     \left\langle q^{- 1} q', p^{- 1} i p i - i p^{- 1} i p \right\rangle = 0\\
     \left\langle q^{- 1} q', p^{- 1} i p + i p^{- 1} i p i \right\rangle = 0
   \end{array}\right. \]
Writing $p = p_0 + i p_1 + j p_2 + k p_3$ we have
\[ u := p^{- 1} i p i - i p^{- 1} i p = 4 ( p_0 p_2 + p_1 p_3 ) j - 4 ( p_1 p_2 - p_0
   p_3 ) k \]
and
\[ v := p^{- 1} i p + i p^{- 1} i p i = -u i
 =  4 ( p_1 p_2 - p_0 p_3 ) j + 4 ( p_0 p_2 + p_1
   p_3 ) k . \]
The two right-hand vectors $u,v$ lie in $\mathrm{Span}(j,k)$ and are either linearly independent
(over $\mathbb{R}$) or both zero. So we have two cases:
\begin{itemizedot}
  \item $u=v=0$, i.e.
  \[ 0 = ( p_0 p_2 + p_1 p_3 )^2 + ( p_1 p_2 - p_0 p_3 )^2 = ( p_0^2 + p_1^2 )
     ( p_2^2 + p_3^2 ) \]
  hence $p$ lies in $\mathrm{Span} ( 1, i )$ or $\mathrm{Span} ( j, k )$ and
  conditions in (\ref{remainingconditions}) hold. Gauging $p$ we may assume
  that $p = 1$ or $p = j$, and the K\"{a}hler angle is then $K = - 1$ or $K = + 1$
  respectively, so that the radius $r$ remains constant (and we may as well
  assume $r = 1$). This case corresponds to Hopf surfaces \cite{P}, i.e. inverse
  images of a curve by the Hopf fibration $\S^3 \rightarrow \S^2$. After a
  possible change of variable in $t$ (replacing $t$ by $t + \varphi ( s )$ for
  some function $\varphi$), we may assume that the curve $s \mapsto e_1 ( s )$
  is Legendrian.

  \item $u,v$ are independent vectors and span $j, k$, thus forcing $q^{-
  1} q'$ to lie in $\mathrm{Span} ( 1, i ) \cap \mathrm{Im} \H$; using gauge action,
  we may assume that $q$ is constant, and up to congruence write $q=1$.
  Reverting to complex coordinates, $e_1 = p = \gamma_1 + \gamma_2 j \simeq (
  \gamma_1, \gamma_2 ) \in \S^3 \subset \mathbb{C}^2$ and $e_2 = p i \simeq ( i \gamma_1, -
  i \gamma_2 )$, while $K = | \gamma_2 |^2 - | \gamma_1 |^2 \neq 0$. We may
  normalize, assuming that $K > 0$ (if $K < 0$ pick the opposite orientation
  on the surface) and set $( \alpha_1, \alpha_2 ) \assign \frac{1}{\sqrt{K}} (
  \gamma_1, \gamma_2 )$, so that
  \[ X = \frac{\sqrt{C}}{\sqrt{K}} ( \gamma_1 ( s ) e^{i t}, \gamma_2 ( s )
     e^{- i t} ) = \sqrt{C} ( \alpha_1 e^{i t}, \alpha_2 e^{- i t} ) \]
  with $| \alpha_1 |^2 - | \alpha_2 |^2 = - 1$, i.e. $( \alpha_1, \alpha_2 )$
  lies in $\H_1^3$, the unit anti-De Sitter space. Again, up to a change in
  variable, we have a Legendrian curve for the indefinite metric in
  $\mathbb{C}^{1, 1}$, i.e. $\left\langle \alpha_1', i \alpha_1 \right\rangle
  - \left\langle \alpha_2', i \alpha_2 \right\rangle = 0$ (see \cite{CLU} or \cite{CC}).
\end{itemizedot}
Summing up, we have proved the following

\begin{theorem}
  A centered cyclic Lagrangian surface may be locally parametrized, up to
  $U(2)$ congruence, by one the following immersions:

\begin{description}
\item[Type I] (complex extensors):
  \[ \begin{array}{cccl}
       X : & I \times \mathbb{R} / 2 \pi \mathbb{Z} & \to & \mathbb{C}^2\\
       & ( s, t ) & \mapsto & r ( s ) e^{i \phi ( s )} ( \cos t, \sin t )
     \end{array} \]
  \item[Type II] (Hopf type):
  \[ \begin{array}{cccl}
       X : & I \times \mathbb{R} / 2 \pi \mathbb{Z} & \to & \mathbb{C}^2\\
       & ( s, t ) & \mapsto & c e^{i t} ( \gamma_1 ( s ), \gamma_2 ( s ) )
     \end{array} \]
  where $\gamma = ( \gamma_1, \gamma_2 )$ is any Legendrian curve of
  $\mathbb{S}^3$ and $c$ is a real constant,

  \item[Type III] (De Sitter type)
  \[ \begin{array}{cccl}
       X : & I \times \mathbb{R} / 2 \pi \mathbb{Z} & \to & \mathbb{C}^2\\
       & ( s, t ) & \mapsto & c \left( \alpha_1 ( s ) e^{i t}, \alpha_2 ( s )
       e^{- i t} \right)
     \end{array} \]
  where $\alpha = ( \alpha_1, \alpha_2 )$ is any Legendrian curve in the unit
  anti-De Sitter space $\H_1^3$ and $c$ is a real constant.
\end{description}
\end{theorem}

\begin{remark} \label{isotropic}
This analysis applies as well if we do not assume that $X$ is an immersion
but only has an isotropic image fibered by circles (we did not use
the immersion hypothesis). So the same conclusion holds and will be used
in the next Section.
\end{remark}

\begin{remark}
Type I surfaces are a particular case of a class of Lagrangian immersions
 which has been first  described in \cite{C1} where they were called \em complex extensors. \em
\end{remark}

\begin{remark}
In the type III case, it may happen that $K$ is identically $\pm 1$. Then we fall back on
type II with a Legendrian curve that actually reduces to a single point. The image of $X$
is therefore a circle (lying in a complex plane).
\end{remark}

\medskip

\section{The general case}

We now consider a surface $\Sigma$ of $\mathbb{R}^4$ which is foliated by
circles. Locally, $\Sigma$ may be parametrized by the following immersion
(identifying as usual $\mathbb{R}^4$ with $\mathbb{C}^2$):
\[ \begin{array}{cccl}
     Y : & I \times \mathbb{R} / 2 \pi \mathbb{Z} & \to & \mathbb{C}^2\\
     & ( s, t ) & \mapsto &  X ( s, t ) + V ( s ) ,
   \end{array} \]
where $V ( s )$ is a $\mathbb{C}^2$-valued function, and
$X ( s, t ) = r ( s ) ( e_1 ( s ) \cos t + e_2 ( s ) \sin t )$ is a centered surface as
in the previous section. Note that we do not assume \emph{a priori}
that $X$ is Lagrangian nor that it is always an immersion.

As we have $Y_t = X_t$ and $Y_s = X_s + V'$, the assumption that $Y$ is
Lagrangian leads to:
\[
\begin{array}{ccl}
0 & = & \langle Y_s, J Y_t \rangle
= \langle X_s, J X_t \rangle + \langle V',   J X_t \rangle \\
& = & \langle X_s, J X_t \rangle + r \cos t \langle V', J e_2
   \rangle - r \sin t \langle V', _1 \rangle
\end{array}   \]
Recall from the previous section that $\langle X_s, J X_t \rangle$ contains
only terms in $\cos 2 t$ and $\sin 2 t$ and a term independent from $t$. Thus
the immersion $Y$ is Lagrangian if and only if: (i) $X$ is cyclic isotropic
(cf Remark~\ref{isotropic}), and (ii) $V'$ belongs to the symplectic
orthogonal of $\mathrm{Span} ( e_1, e_2 )$. Using Theorem 1 we infer:
\begin{itemize}
  \item for type I surfaces, $\mathrm{Span} ( e_1, e_2 ) = \mathrm{Span} ( e^{i
  \phi} ( 1, 0 ), e^{i \phi} ( 0, 1 ) )$ is Lagrangian, so its symplectic
  orthogonal is itself; hence $V' ( s ) = e^{i \phi ( s )} ( W_1 ( s ), W_2 (
  s ) )$ for some real-valued functions $W_1, W_2$;

  \item for type II surfaces, $\mathrm{Span} ( e_1, e_2 )$ is a complex line, so
  its symplectic orthogonal is the same as its Riemannian orthogonal, which is
  \[\mathrm{Span} ( (\bar{\gamma}_{2},-\bar{\gamma}_{1}),
   (i\bar{\gamma}_{2},- i \bar{\gamma}_{1}) );\]
  $V'$ is determined analogously;

  \item for type III surfaces, one can check than a basis of the symplectic
  orthogonal is $( f_1, f_2 )$, where
  $f_{1} = (|\alpha_{2}|^2 , \alpha_{1} \alpha_{2})$,
  $f_2 = (i |\alpha_{2}|^2 , -i \alpha_{1} \alpha_{2})$.
\end{itemize}
So we conclude this section by the

\begin{theorem}
  A cyclic Lagrangian surface may be locally parametrized, up to $U(2)$ congruence,
  by one the following immersions:
  \begin{description}
  \item[Type I]
  \[ \begin{array}{cccl}
       Y : & I \times \mathbb{R} / 2 \pi \mathbb{Z} & \to & \mathbb{C}^2\\
       & ( s, t ) & \mapsto & r ( s ) e^{i \phi ( s )} ( \cos t, \sin t ) +
       \int_{s_0}^s e^{i \phi ( u )} ( W_1 ( u ), W_2 ( u ) ) du
     \end{array} \]
  where $W_1, W_2$ are real valued; in particular when $\phi$ is constant, $Y$ stays
  within the Lagrangian plane
  $\mathrm{Span} ( e^{i \phi} ( 1, 0 ), e^{i \phi} ( 0, 1 ) )$;

  \item[Type II]
  \[ \begin{array}{cccl}
       Y : & I \times \mathbb{R} / 2 \pi \mathbb{Z} & \to & \mathbb{C}^2\\
       & ( s, t ) & \mapsto & c \left( \gamma_1 ( s ) e^{i t}, \gamma_2 ( s )
       e^{i t} \right) + \int_{s_0}^s W ( u )(\bar{\gamma}_{2}(u),-\bar{\gamma}_{1}(u)) du
     \end{array} \]
  where $\gamma = ( \gamma_1, \gamma_2 )$ is any Legendrian curve of
  $\mathbb{S}^3$, $c$ is a real constant and $W$ a complex valued function;
  in particular if $\gamma$ is constant, then up to congruence, we may assume
  $\gamma=(1,0)$ and the immersion becomes $Y(s,t)=(c e^{i t},V_{2}(s))$. Thus
  the immersed surface is a Cartesian product of a circle with a planar curve;

  \item[Type III]
  \[ \begin{array}{cccl}
       Y : & I \times \mathbb{R} / 2 \pi \mathbb{Z} & \to & \mathbb{C}^2\\
       & ( s, t ) & \mapsto & c \left( \alpha_1 ( s ) e^{i t}, \alpha_2 ( s )
       e^{- i t} \right) + \int_{s_0}^s
       (W |\alpha_{2}|^2 , \bar{W} \alpha_{1} \alpha_{2}) d u
     \end{array} \]
  where $\alpha = ( \alpha_1, \alpha_2 )$ is any Legendrian curve in the unit
  anti-De Sitter space $\H_1^3$, $c$ is a real constant and $W$ a complex
  valued function.
  \end{description}
\end{theorem}

\begin{remark} \label{autrebase}
In the type II case, if the curve $(\gamma_1,\gamma_2)$ is in addition regular,
we may assume that it is parametrized by arc length and another
basis of the orthogonal space to $\mathrm{Span} ( e_1, e_2 )$ is
$ \left( (\gamma'_1,\gamma'_2),(i \gamma'_1, i \gamma'_2) \right)$,
so that the immersion may take the alternative form:
$$ Y(s,t)=(\gamma_1(s) e^{i t}, \gamma_2(s) e^{i t}) +
  \int_{s_0}^s W(u)(\gamma'_1(u),\gamma'_2(u)) du.$$
This is a particular case of Lagrangian immersions which have been recently described in \cite{C2}.
This alternative formula will also be useful in the next sections.
\end{remark}


\bigskip

\section{Application to the self-similar equation}
In this section we study the self-similar equation in the case of cyclic Lagrangian surfaces
and prove the following:

\begin{theorem}
 A Lagrangian cyclic surface of $\C^2$ which is a soliton of the mean curvature flow,
i.e. a solution to the self-similar equation
 $$ H + \lambda X^\perp=0 , $$
for some non-vanishing number $\lambda$ is locally congruent to an equivariant example
 described in \cite{A} (in the terminology of the present article, a centered surface of type I)
 or to the Cartesian product $ \S^1(r) \times {\Gamma}$ of some circle $\S^1(1)$
 with a planar self-shrinking curve ${\Gamma}$. Such curves have been studied in detail
 in \cite{AL}.

\end{theorem}

\noindent \textit{Proof.} The proof deals with the three cyclic cases separately:
we first prove that a self-similar surface of type I must be
centered, thus one of the examples of \cite{A}; then we show that there no self-similar surfaces
of type II except the Clifford torus $\S^1 \times \S^1$ (which is also a type I surface)
and the products of curves. Finally we see that there are no self-similar surfaces
of type III at all.

\bigskip

\noindent \textit{Case 1: type I surfaces.}

\medskip

A type I surface is parametrized by an immersion of the form:
 $$\begin{array}{clcl} X: & I \times \R / 2 \pi \Z &  \to &  \C^2 \\
                    & (s,t) & \mapsto &
  r(s) e^{i \phi(s)}(\cos t, \sin t) +
                 \int_{s_0}^s e^{i \phi(u)} (W_1(u),W_2(u)) du ,    \end{array}$$
where $r(s)>0$. Following \cite{ACR}, we shall use the following notations:
$\gamma(s)=r(s) e^{i \phi(s)}$, and, assuming that $\gamma$ is parametrized by
arc length, we shall also denote $\gamma'(s)= e^{i\theta(s)}$.

We start computing the first derivatives of the immersion:
$$ X_s=\gamma' (\cos t, \sin t) + e^{i \phi} (W_1,W_2), \hspace{3em}
 X_t= \gamma (- \sin t , \cos t),$$
 from which we deduce the expression of the induced metric:
\[
E=|X_s|^2= 1 + |W|^2+2 \cos (\theta - \phi) (W_1 \cos t + W_2 \sin t),
\]
\[
F=\<X_s,X_t\>=r \cos (\theta - \phi) (W_2 \cos t - W_1 \sin t),
\qquad G=|X_t|^2= r^2,
\]
and a basis of the normal space to the surface:
$$ N_t= i \gamma (- \sin t , \cos t), \hspace{3em}
 N_s= i \gamma' (\cos t, \sin t) + i e^{i \phi} (W_1,W_2).$$

We now compute the second derivatives of the immersion, in order to calculate the mean
curvature vector:
$$X_{ss}=\gamma'' (\cos t, \sin t) + i \phi' e^{i \phi} (W_1,W_2),$$
$$X_{st}=\gamma' (- \sin t , \cos t), \hspace{3em}
X_{tt}= \gamma (- \cos t , - \sin t).$$
This implies in particular that:
$$ \<X_{ss}, N_t \> = (W_2 \cos t - W_1 \sin t) \sin (\theta - \phi),$$
$$ \<X_{tt}, N_t \>=0, \hspace{3em}
  \<X_{st},N_t\>=r \sin(\theta - \phi).$$

On the other hand, we have:
$$ \<X, N_t \>= a \cos t + b \sin t,$$
where $a := \<i \gamma, \int_{s_0}^s W_1(u) e^{i \phi(u)} du\>$ and
$b:= \<i \gamma, \int_{s_0}^s W_2(u) e^{i \phi(u)}du \>$.

\medskip

We now assume that the immersion $X$ is self-similar, so there exists a non-vanishing real
number $\lambda$ such that:
$$ \< H , N_t \> + \lambda  \<X, N_t \>=0,$$
which is equivalent to
$$ \<X_{ss}, N_t \> G + \<X_{tt},N_t \> E - 2 \<X_{st} , N_t \> F =
 - 2 \lambda (E G - F^2 ) \< X, N_t \>.$$
In the latter expression, the left hand side term is linear in $\cos t$ and $\sin t$
 and the right hand side term is a polynomial of order 3. Linearizing the latter, we easily see
 that the coefficient of $\cos 2t$ is $a W_2 + b W_1$ and the one of
$\sin 2t$ is $a W_1 - b W_2$.
So either $W_1$ and $W_2$ vanish, or $a$ and $b$ vanish. We are going to show that
actually if $a$ and $b$ vanish, then so do $W_1$ and $W_2$.

\medskip

We first write 
$$a= \<i \gamma, \int_{s_0}^s W_1(u) e^{i \phi(u)} \>=
  r \<i e^{i \phi} , \int_{s_0}^s W_1(u) e^{i \phi(u)} du \>$$
$$ =r \left( -\sin \phi \left( \int_{s_0}^s W_1(u) \cos \phi(u) du \right) +
        \cos \phi \left( \int_{s_0}^s W_1(u) \sin \phi(u) du \right) \right)=0 . $$
Thus the derivative of $a/r$ with respect to $s$ must vanish, which yields:
$$  \phi' \left(
  -\cos \phi \left( \int_{s_0}^s W_1(u) \cos \phi(u) du \right)
 - \sin  \phi \left( \int_{s_0}^s W_1(u) \sin \phi(u) du \right)
  \right) $$
$$+ W_1(s)(- \sin \phi \cos \phi + \cos \phi \sin \phi)=0 . $$
Now either $\phi$ is (locally) constant, the curve $\gamma$ is a straight line
passing through the origin and the image of the immersion $X$ is a piece of a plane
(cf Example 2, page 5 of \cite{ACR}), or one can find find points around which
$\phi'$ does not vanish (locally again). In the latter case we get the following linear system:
$$\left\{ \begin{array}{lcc}
-\sin \phi \left( \int_{s_0}^s W_1(u) \cos \phi(u) du \right) +
        \cos \phi \left( \int_{s_0}^s W_1(u) \sin \phi(u) du \right) &=&0\\
 -\cos \phi \left( \int_{s_0}^s W_1(u) \cos \phi(u) du \right)
 - \sin  \phi \left( \int_{s_0}^s W_1(u) \sin \phi(u) du \right)
   & =&0
\end{array} \right. $$
It follows that $\int_{s_0}^s W_1(u) \cos \phi(u) du$ and
$\int_{s_0}^s W_1(u) \sin \phi(u) du$ must vanish except on an isolated set of points of
 $I$, which in turn implies the vanishing of $W_1$. Analogously it can be shown that
 $W_2$ vanishes as well, so finally the surface must be centered. Lagrangian, self-similar,
 centered surfaces have been described in detail in \cite{A}. We only recall here that such
 surfaces are obtained from planar curves $\gamma$ which are solutions of the following
 equation:
$$ k=\< \gamma, N\>\left(\frac{1}{|\gamma|^2}-\lambda
 \right),   $$
where $k$ is the curvature of $\gamma$ and $N$ its unit
normal vector. This equation admits a countable family of closed solutions
    which is parametrised by two relatively prime numbers $p$
    and $q$ subject to the condition
    $ p/q \in (1/4, 1/2);$
    $p$ is the winding number of the curve and $q$ is the number
    of maxima of its curvature. Except for the circles, none of these curves is embedded.
\bigskip

\noindent \textit{Case 2: type II surfaces.}

\medskip
As our discussion is local, we consider the two following cases: either the Legendrian
curve $\gamma$ is regular, or it reduces to a single point, and then the surface is
a product of a circle with some plane curve. In the first case, the surface
may parametrized by an immersion of the form (cf Remark~\ref{autrebase}):
$$\begin{array}{clcl} X: & I \times \R / 2 \pi \Z &  \to &  \C^2 \\
                    & (s,t) & \mapsto &
 c(\gamma_1(s) e^{i t}, \gamma_2(s) e^{i t}) +
  \int_{s_0}^s W(u)(\gamma'_1,\gamma'_2) du
      \end{array}$$
 where $(\gamma_1(s),\gamma_2(s))$ is some unit speed Legendrian curve of $\S^3$.
Without loss of generality we fix $c=1$.
We start computing the first derivatives of the immersion:
$$ X_s=(\gamma'_1 e^{i t}, \gamma'_2 e^{i t}) + W(\gamma'_1,\gamma'_2), \hspace{3em}
 X_t=(i \gamma_1 e^{i t}, i \gamma_2 e^{i t}),$$
 from which we deduce the expression of the induced metric:
 $$E=|X_s|^2= 1 + |W|^2 + 2 \<W, e^{i t} \>, \qquad
F=\<X_s,X_t\>=0, \qquad G=|X_t|^2=1,$$
and a basis of the normal space to the surface:
$$N_s=(i \gamma'_1 e^{i t}, i \gamma'_2 e^{i t}) + W(i\gamma'_1,i\gamma'_2), \hspace{3em}
 N_t=(- \gamma_1 e^{i t}, - \gamma_2 e^{i t}) .$$

We now compute the second derivatives of the immersion, in order to calculate the mean
curvature vector:

$$X_{ss}=(\gamma''_1 e^{i t}, \gamma''_2 e^{i t}) + W(\gamma_1'', \gamma''_2)
  +  W'(\gamma'_1, \gamma'_2),$$
$$ X_{st}=(i\gamma'_1 e^{i t}, i \gamma'_2 e^{i t}), \hspace{4em}
X_{tt}=(-\gamma_1 e^{i t}, - \gamma_2 e^{i t}). $$

We first notice that
 $$ \<X, N_t \> =
    -1 + a \cos t + b \sin t,$$
where $a$ and $b$ depend only on the variable $s$.

On the other hand,
$$ 2 \<H, N_t \> = \frac{ \<X_{ss},N_t \>}{E}+  \frac{ \<X_{tt},N_t \>}{G}
 = \frac{ \<X_{st},N_s \>}{E}+  \frac{ \<X_{tt},N_t \>}{G}$$
$$ =\frac{1 + \<W, e^{i t} \>}{1 + |W|^2 + 2 \<W, e^{i t} \>} +1
      = \frac{2+ |W|^2 +3\<W, e^{i t} \>}{1 + |W|^2 + 2 \<W, e^{i t} \>}.$$
Thus the equation
$\<H, N_t \> + \lambda \<X, N_t \> = 0$ holds
for some non-vanishing constant $\lambda$ if and only if $a$ and $b$ vanish and $|W|=1$
or $W=0$. In particular, $\<X,N_t \> = -1$.

We leave to the reader the easy task to check that if $W$ vanishes, the immersion is
self-similar if and only if $\<\gamma'', J \gamma'\>$ vanishes, that is the curve $\gamma$ has
vanishing curvature and thus is a great circle. The corresponding Lagrangian surface is
the Clifford torus $\frac{1}{\sqrt{2}} \S^1 \times \frac{1}{\sqrt{2}} \S^1$.
So we assume
in the remainder that $|W|=1$. In particular, we may write $W=e^{i\phi}$ and
$W'= i \phi' e^{i \phi}$, where $\phi$ is some real function of the variable $s$.

\bigskip

We want to look at the other scalar equation $\<H, N_s \> + \lambda \<X, N_s \> = 0$,
so we compute
$$ \<X_{ss},N_s \>
 =\< \gamma'', J\gamma' \> ( 1 + |W|^2) + \<i W, W'\> + 2 \<W, e^{-i t}\> \< \gamma'', J\gamma' \>
 + \<W', e^{i t}\>,$$
$$ \<X_{tt}, N_s\> =0.$$
Hence we get
\[ \begin{array}{rcl}
2 \<H, N_s \> &=& \displaystyle
\frac{ \<X_{ss},N_s \>}{E}+  \frac{ \<X_{tt},N_s \>}{G}
 = \frac{\<X_{ss},N_s \>}{1 + |W|^2 + 2 \<W, e^{i t} \>} + 0
 \\
& =& \displaystyle
\frac{\< \gamma'', J\gamma' \> ( 1 + |W|^2) + \<i W, W'\> +2 \<W, e^{-i t}\> \< \gamma'',J\gamma' \>
 + \<W', e^{i t}\>}{1 + |W|^2 + 2 \<W, e^{i t} \>}
 \\
&=&  \displaystyle \frac{ 2 \<\gamma'', J\gamma' \>  + \phi' +
  2 \cos (\phi+ t)\< \gamma'',J\gamma' \> - \phi' \sin(\phi - t) }{2  + 2 \cos (\phi - t)}
  \\
&=& \displaystyle
\frac{2 \<\gamma'', J\gamma' \>  + \phi' + (2 \< \gamma'',J\gamma' \> \cos \phi - \phi' \sin \phi) \cos t}
{2 (1 + \cos \phi \cos t + \sin \phi \sin t)}
\\
&& + \displaystyle
 \frac{ (- 2 \< \gamma'',J\gamma' \> \sin \phi + \phi' \cos \phi) \sin t}{
 2 (1 + \cos \phi \cos t + \sin \phi \sin t)}
 \end{array} \]

On the other hand, it is easy to see that $\<X,N_s \>$ takes the form
$a \cos t + b \sin t + c$, where $a$ $b$ and $c$ depend on the variable $s$.
 This forces  $\<H, N_s \>$ to take a simpler form and implies
that
 $$ \frac{  2 \<\gamma'', J\gamma' \>  + \phi'}{2}=
 \frac{2 \< \gamma'',J\gamma' \> \cos \phi - \phi' \sin \phi}{2 \cos \phi}
= \frac{- 2 \< \gamma'',J\gamma' \> \sin \phi + \phi' \cos \phi}{ 2\sin \phi}$$
In particular we have $\phi'= - \phi' \tan \phi$, which implies that $\phi$ is constant.
It follows that $W$ is constant as well, so $\< \gamma'', J\gamma' \>$ vanishes and the curve $\gamma$
has vanishing curvature, so it is a great circle of the unit sphere. There is no loss of
generality to assume that
$\gamma(s)=\frac{1}{\sqrt{2}}(e^{is}, e^{-is})$.
 Now the immersion takes the following, explicit form:
 $$ X(s,t) = \frac{1}{\sqrt{2}}(e^{i(s+t)},e^{i(-s+t)}) +
                \frac{1}{\sqrt{2}} (-i e^{i(\phi+s)},i e^{i (\phi-s)}).$$
But in this case it is easy to check $\<X,N_t \>=-1$ does not hold, so we conclude that there
is no self-similar type II surface with a regular curve $\gamma$ and non-vanishing $W$.

\bigskip
It remains to treat the case of a Cartesian product $\S^1(r)
\times  \Gamma$ of a circle  with an arbitrary curve
$\Gamma$. It is straightforward that such a product is
self-similar if and only if both curves are solutions of the
equation
 $$ k + \lambda \< X , N \> =0$$
 for the same $\lambda$.
A circle of radius $r$
is trivially self-similar for $\lambda = r^{-2}$.
The other self-similar
 curves with positive $\lambda$ (self-shrinking curves)
have been described in \cite{AL}. Except for the circles, none of
them is embedded.

\bigskip

\noindent \textit{Case 3: type III surfaces.}

\medskip

A type III surface is parametrized by an immersion of the form:
$$\begin{array}{clcl} X: & I \times \R / 2 \pi \Z &  \to &  \R^4 \\
                    & (s,t) & \mapsto &
 c(\alpha_1(s) e^{i t}, \alpha_2(s) e^{-i t})+
      \int_{s_0}^s (W|\alpha_2|^2,\bar{W} \alpha_1 \alpha_2)du,
      \end{array}$$
where $(\alpha_1,\alpha_2)$ is a Legendrian curve of $\H_1^3$, in particular we may assume
that
 $$ |\alpha'_1|^2 - |\alpha'_2|^2 = 1  \quad
        \<\alpha'_1, i \alpha\> - \< \alpha'_2,i \alpha_2 \>=0.$$
Observe that this implies the following identity (cf \cite{CC}): $|\alpha_1|=|\alpha'_2|$.
Again, without loss of generality we fix $c=1$.
\medskip

We start computing the first derivatives of the immersion:
$$ X_s=(\alpha'_1 e^{i t}, \alpha'_2 e^{-i t}) + (W|\alpha_2|^2,\bar{W}\alpha_1 \alpha_2),
\quad X_t=(i \alpha_1 e^{i t}, -i \alpha_2 e^{-i t}),
$$
from which we deduce the expression of the induced metric: 
 $$E = 1 + 2 |\alpha_1|^2 + |W|^2 (1+ |\alpha_1|^2)(1+ 2|\alpha_1|^2)
 +  |\alpha_2|^2 \< \alpha'_1 \bar{\alpha_1},  W  e^{-i t} \>
        + |\alpha_1|^2 \<  \alpha'_2 \bar{\alpha_2}, \bar{W} e^{i t}\>,  $$
 $$ F= -2  (1 + |\alpha_1|^2)( \Im( \alpha_1 \bar{W})\cos t
                                + \Re (\alpha_1 \bar{W})\sin t),$$
$$ G=1 + 2|\alpha_1|^2, $$
and a basis of the normal space to the surface:
$$ N_s=(i\alpha'_1 e^{i t}, i\alpha'_2 e^{-i t}) +
                    (i W|\alpha_2|^2,i\bar{W}\alpha_1 \alpha_2),$$
$$ N_t=(- \alpha_1 e^{i t}, \alpha_2 e^{-i t}).$$
 We now compute the second derivatives of the immersion:
$$X_{ss}=(\alpha''_1 e^{i t}, \alpha''_2 e^{-i t}) +
       \frac{d}{d s}(W|\alpha_2|^2,\bar{W} \alpha_1 \alpha_2),$$
$$X_{st}= (i\alpha'_1 e^{i t}, -i\alpha'_2 e^{-i t}), \hspace{3em}
 X_{tt}= (- \alpha_1 e^{i t},- \alpha_2 e^{-i t});$$
and some coefficients of the second fundamental form:
$$ \<X_{tt}, N_t \>
= |\alpha_1|^2 - |\alpha_2|^2
=-1,  \hspace{3em} \<X_{st}, N_t \>
 = 2 \<i\alpha'_1, \alpha_1 \>; $$
 Moreover it is easy to see that the coefficient
 $\<X_{ss}, N_t \>$ is an affine function of $\cos t$ and
 $\sin t$.

 \medskip

On the other hand, we have
$$ \<X, N_t \>    = 1+ a \cos t + b \sin t,$$
where $a$ and $b$ depend only on the variable $s$.

\medskip

 We now assume that
$\<H, N_t \> + \lambda \<X, N_t \> = 0$
holds for some non-vanishing constant $\lambda$; this implies
$$ G\<X_{ss},N_t \>+ E \<X_{tt},N_t \>
     - 2 F \<X_{st}, N_t\> =
  -\lambda (1+ a \cos t + b \sin t)(E G - F^2).$$
In the latter expression, the left hand side term is linear in
$\cos t$ and $\sin t$ and the right hand side is a polynomial of
degree 3. Linearizing the latter and introducing the notation $F=
\tilde{a} \cos t + \tilde{b} \sin t$, we easily see that the
coefficient of $\cos 3t$ is (up to a multiplicative constant) $a
\tilde{a}^2 - 2b \tilde{a}\tilde{b}$ and the one of $\sin 3t$ is $
2a\tilde{a}b-b \tilde{b}^2$. Thus we get the following system:
$$ \left\{ \begin{array}{ccc}    a \tilde{a}^2 - 2b \tilde{a}\tilde{b} & = &0 \\
  2a\tilde{a}b- b \tilde{b}^2 & = &0
 \end{array} \right. $$
This implies that either $\tilde{a}$ and $\tilde{b}$ vanish, or $a$ and $b$ vanish. However, in the latter case, the
coefficient of $\cos 2t$ is (up to a multiplicative constant)
$\tilde{a}^2 + \tilde{b}^2$ and the one of $\sin 2t$ is $2 \tilde{a} \tilde{b}$ thus again $\tilde{a}$
and $\tilde{b}$ must vanish. We conclude by observing that
 $$\tilde{b} +i \tilde{a}=-2 (1 + |\alpha_1|^2) \alpha_1 \bar{W},$$
 so  either $\alpha_1$ or $W$ vanishes.
In the first case, as the discussion is local we may assume
that $\alpha_1$ vanishes identically, and the immersion takes the form
$$ Y(s,t) =(V_1(s), \alpha_2 e^{-i t}),$$
so the immersed surface is a product of curves; this case was already treated
in the former section (type II surfaces).

\medskip

To complete the proof it remains to show that a centered type III surface cannot be
self-similar. We shall use the following result which can be found in \cite{CC}
(Proposition 2.1, page 3 and Corollary 3.5, page 9) and that we state here in accordance to our
own notations:

\begin{proposition} (\cite{CC})
Let $\gamma=(\gamma_1,\gamma_2)$ be unit speed Legendrian curve in $\S^3$ and $\alpha=(\alpha_1,\alpha_2)$
a unit speed Legendrian curve in $\H_1^3$.

Then the following immersion
$$\begin{array}{clcl} X: & I_1 \times I_2 &  \to &  \C^2 \\
                    & (s,t) & \mapsto &
 (\alpha_1(s) \gamma_1(t), \alpha_2(s) \gamma_2(t)),
      \end{array}$$
is conformal and Lagrangian; moreover, its mean curvature vector is given by
$$ H = e^{-2u}(k_{\alpha} J X_s + k_{\gamma} J X_t),$$
where $e^{-2u}$ is the conformal factor and $k_{\alpha}$ and $k_{\gamma}$ are the curvature
functions of $\alpha$ and $\gamma$ respectively.
\end{proposition}

By taking $\gamma(t)=(\frac{1}{\sqrt{2}} e^{i t}, \frac{1}{\sqrt{2}} e^{-i t})$
in the immersion above, we recover exactly a centered type III immersion.
In particular, $\gamma$ is a great
circle of $\S^3$ and has vanishing curvature. This implies that $\<H, J X_t \>$ vanishes.
On the other hand, we calculate
$$ \< X, J X_t \> = \frac{1}{2}
  \< (\alpha_1 e^{i t}, \alpha_2 e^{-i t}),(-\alpha_1 e^{i t}, \alpha_2 e^{-i t}) \>
 = \frac{1}{2} (-|\alpha_1|^2 + |\alpha_2|^2) = \frac{1}{2} .$$
So we deduce that this immersion $X$ cannot be solution of the self-similar equation
$H + \lambda X^\perp=0$.

\bigskip

\section{Application to the Hamiltonian stationary equation}
 In this section, we study the Hamiltonian stationary equation. We first recall
  that to a Lagrangian surface $\Sigma$ is attached
 the \emph{Lagrangian angle function} which is
defined by the formula
 $$ {\det}_{\C} (e_1,e_2) = e^{i \beta},$$
where $(e_1,e_2)$ is a (local) orthonormal tangent frame $\Sigma$
(and thus a Hermitian frame of $\C^2$). Next the surface $\Sigma$
is said to be \emph{Hamiltonian stationary}
 if it is a critical point of the area with respect to Hamiltonian
 variations, i.e. variations generated by vector fields $V$ such
 that $V \lrcorner \omega$ is exact. The Euler--Lagrange
 equation of this variational problem is
 $$ \Delta \beta=0.$$
In other words, a Lagrangian surface is Hamiltonian stationary
 if and only if its Lagrangian angle is harmonic with respect to
 the induced metric.

 In this section, we shall characterize
 Hamiltonian stationary type I and type II surfaces. The situation
 is somewhat similar to the self-similar case treated in the previous
 section: there are examples of non-trivial Hamiltonian stationary
 centered type I surfaces,
 which are described in Subsection 4.1; next we show in Subsection
 4.2 that there are no Hamiltonian stationary type I surfaces
 which are not centered. Finally, we prove in Subsection 4.3 that
 the only Hamiltonian stationary type II surfaces are
 Cartesian products of circles. Similar calculations can be pursued
  in the type III case, however the dependency of the
translation term on the Legendrian curve $\alpha$ lying in
$\H^3_1$ makes explicit computations rapidly too complicated. 
We conjecture nonetheless that Hamiltonian stationary type III surfaces
are centered. Finally we point out that those have been described 
explicitly in \cite{CC}, Section 4.3, page 13: a centered type III surface 
is Hamiltonian stationary if and only if the curvature of the generating 
curve $\alpha$ is an affine function of its arclength parameter.

\subsection{The centered type I case }
We consider an immersion of a centered type I surface, as in Section 1:
\[ \begin{array}{cccl}
       X : & I \times \mathbb{R} / 2 \pi \mathbb{Z} & \to & \mathbb{C}^2\\
       & ( s, t ) & \mapsto & r ( s ) e^{i \phi ( s )} ( \cos t, \sin t )
     \end{array} \]
where the planar curve $\gamma(s)= r ( s ) e^{i \phi ( s )}$ is assumed to be parametrized by arc length.
It follows that there exists $\theta(s)$ such that $\gamma'(s)= e^{i \theta(s)}.$ Introducing
the variable $\alpha= \theta - \phi$,
a straightforward computation (cf \cite{ACR}) yields that the induced metric of the immersion is
$$ g = \left(\begin{array}{cc} 1 & 0 \\ 0 & r^2 \end{array} \right)$$
and its Lagrangian angle function is
 $\beta= \theta + \phi= \alpha + 2 \phi$;
it follows
that
$$ \Delta \beta = 0 \Leftrightarrow \partial_s ( g^{11} \sqrt{\det g} \beta_s)= 0$$
$$ \Leftrightarrow r (\alpha' + 2 \phi')= C,$$
where $C$ is some real constant.
So we are left with the following differential system
$$ \left\{ \begin{array}{ccc}
 r'&=& \cos \alpha \\ \alpha'&=& \frac{C- 2 \sin \alpha}{r} \end{array} \right.$$
which admits a first integral:
$$ E(r,\alpha)= r^2 (C-2\sin \alpha )$$

\begin{figure}  \begin{center}
\includegraphics[width=5cm]{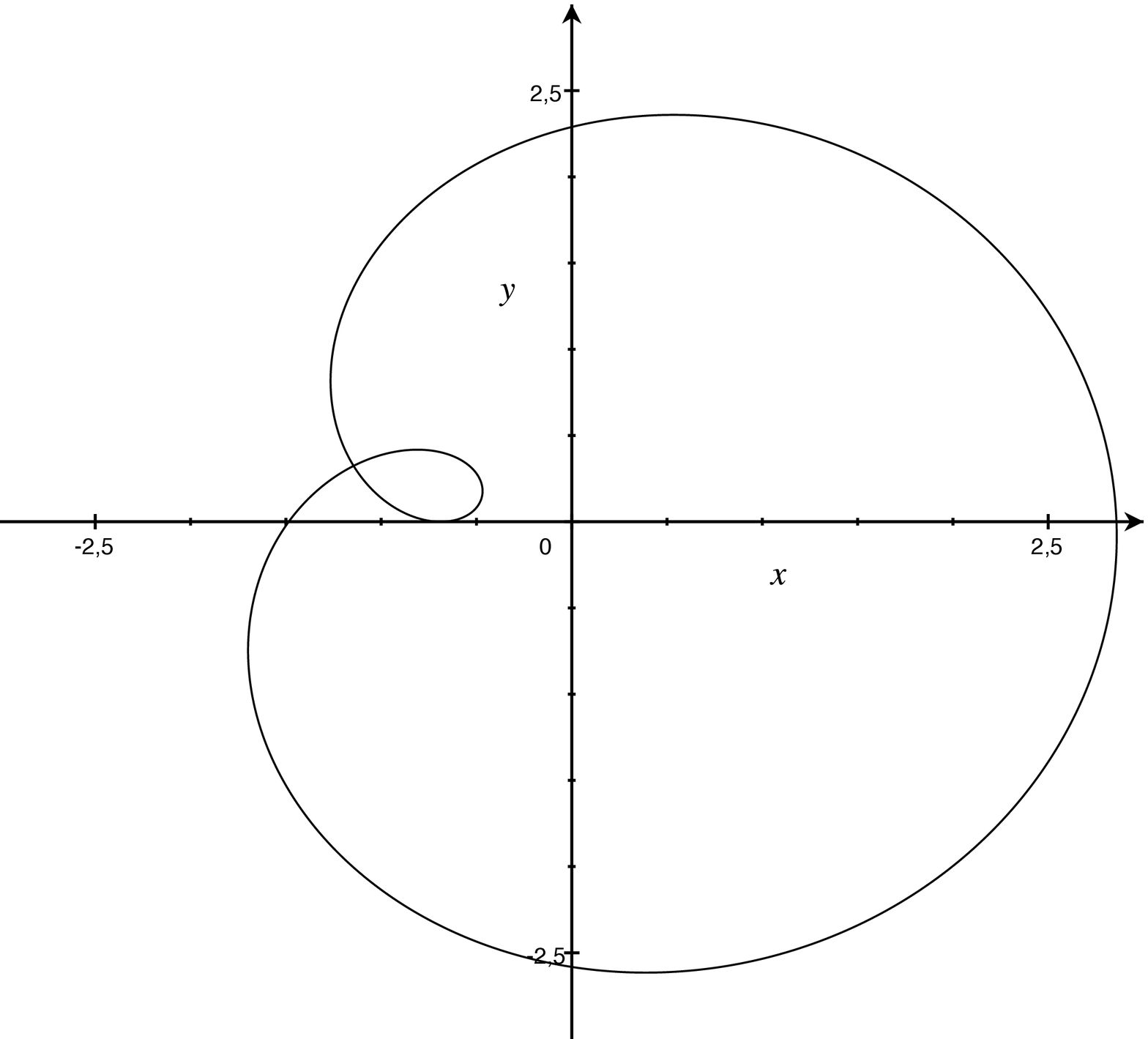} 
\includegraphics[width=5cm]{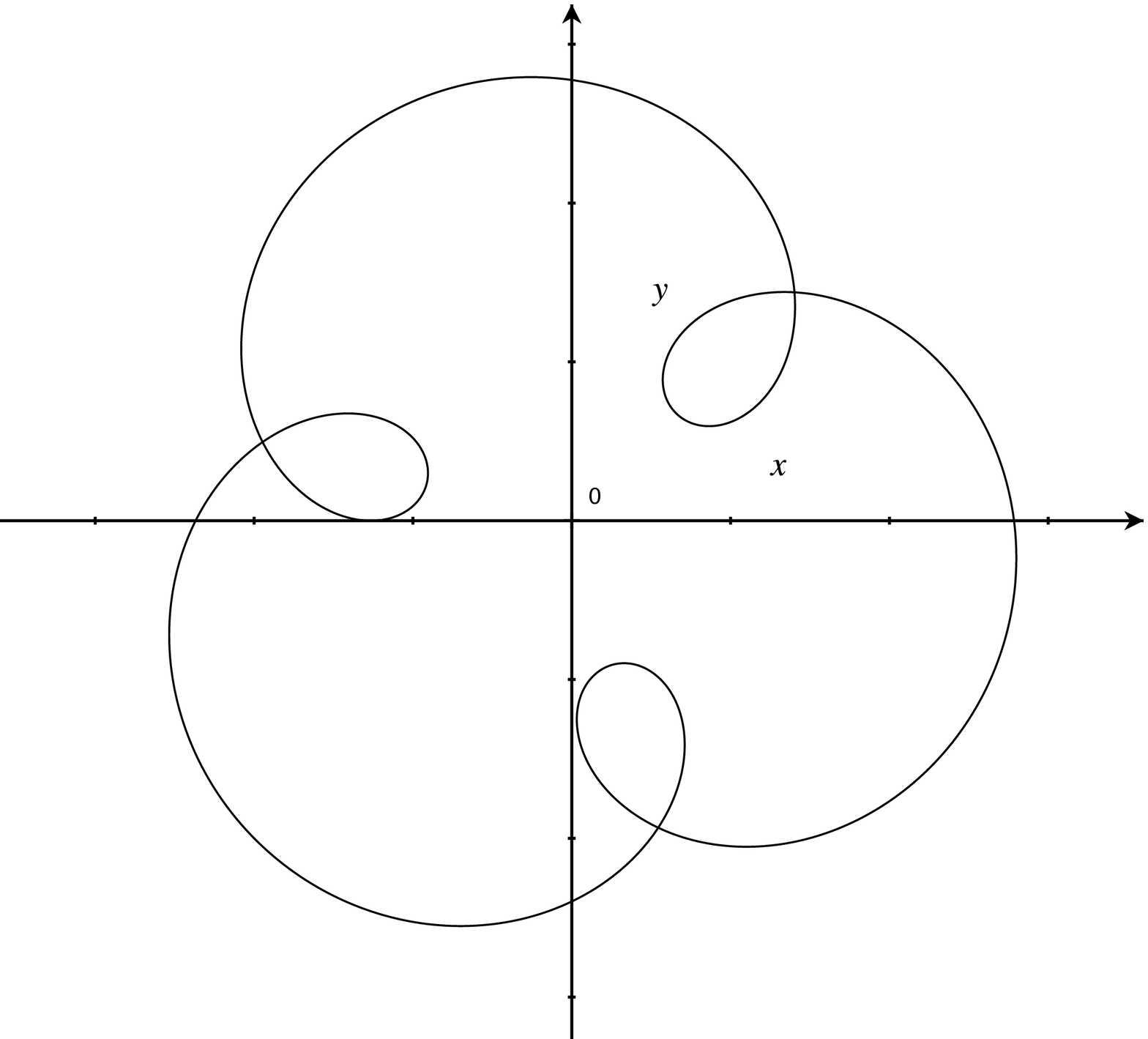} \\
\includegraphics[width=5cm]{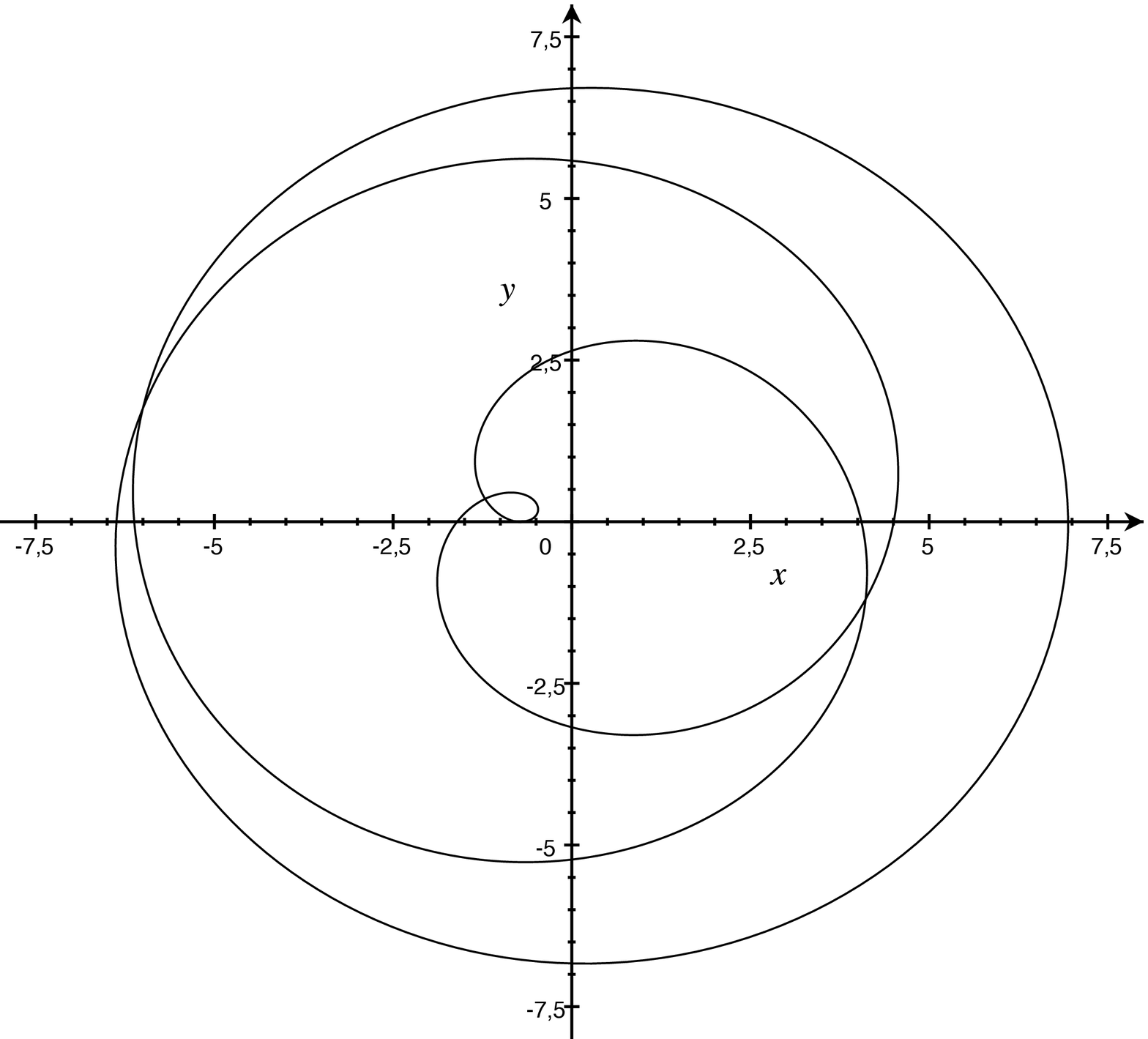}
   \caption{Closed curves for $C>2$ with respective total angular variation $\Phi=2 \pi, 
   2 \pi/3, 6 \pi$.}
\end{center}
\end{figure}
\begin{figure}  \begin{center}
\includegraphics[width=14cm, height=10cm]{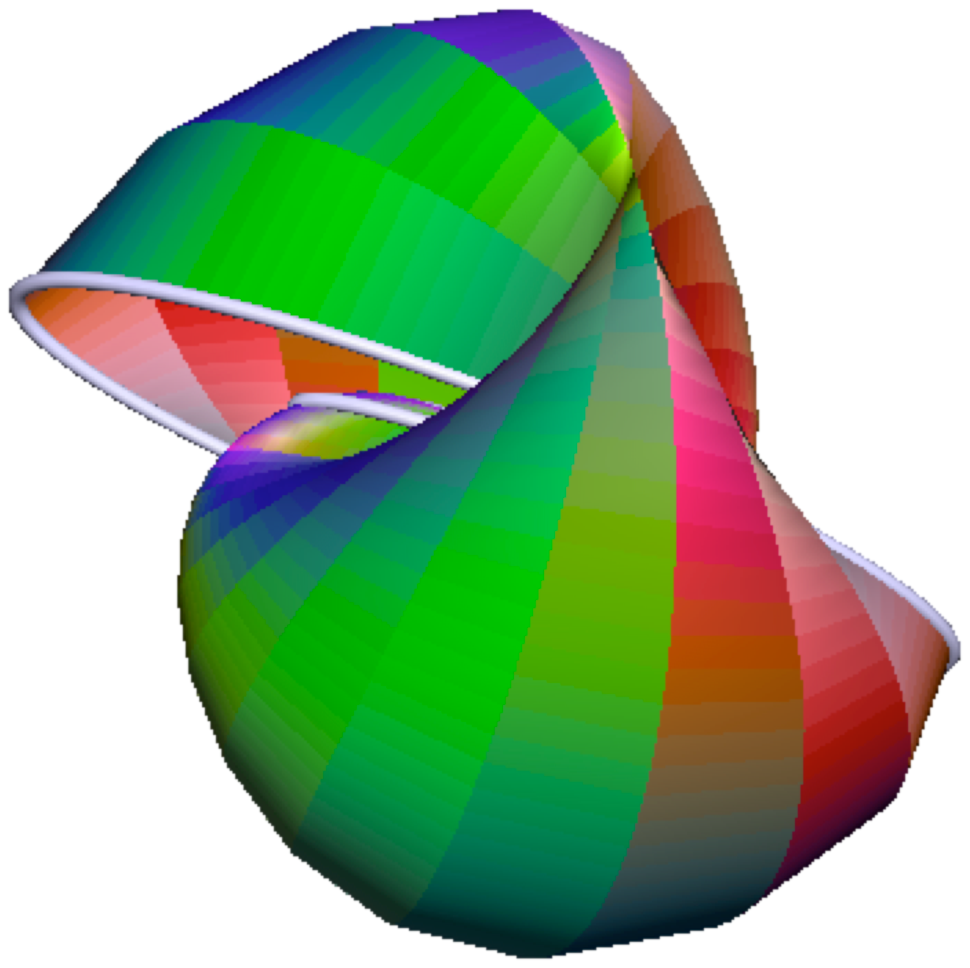} 
   \caption{Surface corresponding to the first closed curve, partially cut to show the circles in grey, 
   projected down to $\R^3$.}
\end{center}
\end{figure}

\noindent \textbf{First case: $|C|>2$.}
Here the trajectories are bounded in the variable $r$ since
$$r = \sqrt{\frac{E}{C-2 \sin \alpha }}.$$
therefore the curves $\gamma$ are also bounded. In order to know whether
such curves are closed, we calculate the variation of the angle of the curve with respect
to the origin along a period of the variable $\alpha$:
\[ \Phi(C) = \int \frac{\sin \alpha}{r} d s = \int_0^{2 \pi} \frac{\sin \alpha}{C - 2 \sin
\alpha}d \alpha . \] 
An easy computation shows that $\lim_{C \to \pm
2}=+\infty$ and $\lim_{C \to \infty}=0$. It follows that for all rational numbers 
$p/q$ there exists  $C(p/q)$ such that the corresponding curve
$\Phi(C_q)= 2 p \pi/ q$. By repeating $q$ times this pattern, we obtain 
a $q$-symmetric curve. It is easy to check that these curves are never embedded.

\bigskip

\noindent \textbf{Second case: $|C|<2$.}
Here the trajectories are unbounded. However we still have
$$ \Phi(C) = \int \frac{\sin \alpha}{r} d s = \int_{\alpha_-}^{\alpha_+} \frac{\sin \alpha}{C - 2 \sin \alpha} d \alpha,$$
and
$$  \Phi(C) = \int_{\alpha_+}^{\alpha_- + 2 \pi} \frac{\sin \alpha}{C - 2 \sin \alpha}d \alpha,$$
where $\alpha_{-} < \alpha_{+}$ are defined to be the two roots of $ 2
\sin \alpha_{\pm} = C$. In both cases, if $C$ does not vanish,
 $\Phi(C)= \infty$, therefore the curve $\gamma$ has two ends which spiral. The case of vanishing $C$
implies $\beta$ to be constant, so this is the special Lagrangian
case, which have already been treated in \cite{A}, Section 4 (cf
also \cite{CU}). It is proved that we get, up to congruences, a
unique surface, the \emph{Lagrangian catenoid}.



\bigskip

\noindent \textbf{Third case: $|C|=2$.}

Here all the points of the vertical lines  $\alpha= 0 \bmod \pi$
are equilibrium points. The corresponding curves $\gamma$ are round circles centered
at the origin. The corresponding surfaces are products of circles $\S^{1}(r) \times \S^{1} (r)$,
sometimes called Clifford tori. In the region $\alpha \neq 0 \bmod \pi$,
the situation is analogous the second case: the curves $\gamma$ have two spiraling ends.

\subsection{Characterization of Hamiltonian stationary Type I surfaces}

\begin{proposition}
 A non-planar Hamiltonian stationary Lagrangian cyclic surface of $\C^2$ of type I
  must be centered, and thus is one of the surfaces described in Section~4.1
\end{proposition}

\noindent \textit{Proof.} The proof is done by contradiction. We consider
a type I immersion
\[ \begin{array}{cccl}
       Y : & I \times \mathbb{R} / 2 \pi \mathbb{Z} & \to & \mathbb{C}^2\\
       & ( s, t ) & \mapsto & r ( s ) e^{i \phi ( s )} ( \cos t, \sin t ) +
       \int_{s_0}^s e^{i \phi ( u )} ( W_1 ( u ), W_2 ( u ) ) du
     \end{array} \]
with non-vanishing  $W=(W_1,W_2)$ and we show that the Hamiltonian
stationarity leads to a contradiction. We recall some notations
already introduced  in \cite{ACR}:
\[ x = (\cos t, \sin t), \]
\[ A \assign \left( 1 + \left\langle W, x \right\rangle^2 + 2 \cos \alpha\langle W, x \rangle \right)^{- 1 / 2} . \]
and
\[ B := k + \left( k \cos \alpha + \frac{\sin \alpha \cos \alpha}{r}
\right)  \left\langle W, x \right\rangle - \sin \alpha
\left\langle W', x \right\rangle + \frac{\sin \alpha}{r}
\left\langle W, x \right\rangle^2.\] 
We also exploit the
computations done in \cite{ACR}, Section 4, where it was shown
that
$$ A^{-6} \Delta \beta = \mathrm{I} + \mathrm{II} + \mathrm{III} + \mathrm{IV} 
+ \mathrm{V} + \mathrm{VI} +\mathrm{VII}$$
where
\begin{eqnarray*}
  \mathrm{I} & := &3 B \left( \cos \alpha
  \left\langle W', x \right\rangle - \alpha' \sin \alpha
  \left\langle W, x \right\rangle + \left\langle W, x \right\rangle
  \left\langle W', x \right\rangle \right)
\\
  \mathrm{II} & := & - A^{- 2} \left[ B' - \frac{\sin
  \alpha}{r}  \left( \cos \alpha \left\langle W', x \right\rangle -
  \alpha' \sin \alpha \left\langle W, x \right\rangle + \left\langle W, x
  \right\rangle  \left\langle W', x \right\rangle \right) \right]
\\
  \mathrm{III} &:=  & - A^{- 4} \partial_s \left(
  \frac{\sin \alpha}{r} \right)
\\
  \mathrm{IV} & := & - 3 B \frac{\cos \alpha +
  \left\langle W, x \right\rangle}{r} (|W|^2 - \left\langle W, x
  \right\rangle^2)
\\
  \mathrm{V} & := &  \frac{A^{- 2}}{r} \bigg[ \left( k \cos
  \alpha + \frac{\sin \alpha}{r}  \left\langle W, x \right\rangle \right)
  (|W|^2 - \left\langle W, x \right\rangle^2)  \\
	&&+ \sin \alpha \left(
  \left\langle W', W \right\rangle - \left\langle W', x
  \right\rangle  \left\langle W, x \right\rangle \right) \bigg]
\\
  \mathrm{VI} & := & - A^{- 2}  \frac{\sin \alpha}{r^2}
  (\cos \alpha + \left\langle W, x \right\rangle) \left( |W|^2 - \left\langle
  W, x \right\rangle^2 \right) - 2 A^{- 4}  \frac{\sin \alpha}{r^2}
  \left\langle W, x \right\rangle
\\
  \mathrm{VII} & := & - A^{- 2}  \frac{B}{r} (\cos \alpha
  + \left\langle W, x \right\rangle) - A^{- 4}  \frac{\sin \alpha \cos
  \alpha}{r^2}
\end{eqnarray*}
 For
fixed $s$ we may view $A^{- 6} \Delta \beta$ as a polynomial in
$x$; then the highest degree term is not anymore in $\left\langle
W, x \right\rangle^5$ as in \cite{ACR} but only in the fourth
power. Forgetting all powers of $x$ beneath the fourth, we have:
\begin{eqnarray*}
  \mathrm{I} & \equiv & \frac{3 \sin \alpha}{r}  \left\langle W', x
  \right\rangle  \left\langle W, x \right\rangle^3 
\\
  \mathrm{II} & \equiv & - \partial_s \left( \frac{\sin \alpha}{r} \right)  \left\langle W, x
  \right\rangle^4 - \frac{\sin \alpha}{r}  \left\langle W', x
  \right\rangle   \left\langle W, x \right\rangle^3
\\
  \mathrm{III} & \equiv &- \partial_s \left( \frac{\sin \alpha}{r} \right)  \left\langle W, x
  \right\rangle^4
\\
 \mathrm{IV} & \equiv & \frac{3 \sin \alpha}{r^2}  \left\langle W, x \right\rangle^5 + \left(
  \frac{3 k \cos \alpha}{r} + \frac{6 \sin \alpha \cos \alpha}{r^2} \right)
  \left\langle W, x \right\rangle^4 \\
	&& - \frac{3 \sin \alpha}{r} \left\langle
  W', x \right\rangle  \left\langle W, x \right\rangle^3 
\\
  \mathrm{V} & \equiv & - \frac{\sin \alpha}{r^2}  \left\langle W, x \right\rangle^5 -
  \left( \frac{k \cos \alpha}{r} + \frac{2 \sin \alpha \cos \alpha}{r^2}
  \right)  \left\langle W, x \right\rangle^4 \\
	&& - \frac{\sin \alpha}{r}
  \left\langle W', x \right\rangle  \left\langle W, x \right\rangle^3
\\
  \mathrm{VI} & \equiv & - \frac{\sin \alpha}{r^2}  \left\langle W, x \right\rangle^5 -
  \frac{5 \sin \alpha \cos \alpha}{r^2}  \left\langle W, x \right\rangle^4 
\\
 \mathrm{VII} & \equiv & - \frac{\sin \alpha}{r^2}  \left\langle W, x \right\rangle^5 - \left(
  \frac{5 \sin \alpha \cos \alpha}{r^2} + \frac{k \cos \alpha}{r} \right)
  \left\langle W, x \right\rangle^4 \\
	&& + \frac{\sin \alpha}{r} \left\langle
  W', x \right\rangle  \left\langle W, x \right\rangle^3
\end{eqnarray*}

At given $s$, for $\Delta \beta$ to be zero for all $t$, we need (replacing
$k$ by $\alpha' + \sin \alpha / r$ and using $r' = \cos \alpha$,
$\phi' = \sin \alpha / r$),
\begin{eqnarray*}
  0 & = & \frac{\sin \alpha}{r}  \left\langle W', x \right\rangle +
  \left( 2 \partial_s \left( \frac{\sin \alpha}{r} \right) + \frac{5 \sin
  \alpha \cos \alpha}{r^2} - \frac{\alpha' \cos \alpha}{r} \right)
  \left\langle W, x \right\rangle\\
  & = & \frac{\sin \alpha}{r}  \left\langle W', x \right\rangle + \left(
  \partial_s \left( \frac{\sin \alpha}{r} \right) + \frac{4 \sin \alpha \cos
  \alpha}{r^2} \right)  \left\langle W, x \right\rangle
\end{eqnarray*}
for almost all values, hence for all. If $W$ and $W'$ are not colinear,
then $\left\langle W, x \right\rangle$ and $\left\langle W', x
\right\rangle$ are linearly independent as functions of $t$, hence
\[ \frac{\sin \alpha}{r} = \partial_s \left( \frac{\sin \alpha}{r} \right) +
   \frac{4 \sin \alpha \cos \alpha}{r^2} = 0 \]
so that $\theta \equiv \phi \bmod  \pi$, $\phi$ is constant
(and it is half the Lagrangian angle). The solution is then
trivial: a subset of some Lagrangian plane, spanned by circles.

If $W'$ is colinear to $W$ for an open interval in the $s$
variable, then $W (s) = w (s) W^0$ for some real valued function
$w (s)$, and the center of the circle at $s$ lies in the plane $\C
W^0$. Using real rotations, and up to a reparametrization in $t$,
we may as well assume that $W^0 = (1, 0)$ and we arrive at the
following equation
\[ \frac{\sin \alpha}{r}  w' + \left( \partial_s \left( \frac{\sin
   \alpha}{r} \right) + \frac{4 \sin \alpha \cos \alpha}{r^2} \right) w = 0 .
\]
whose solution (up to scaling) is $w = c (r^3 \sin \alpha)^{- 1}, $
where $c$ is some non vanishing constant (the case of vanishing $c$ would correspond to the centered case).
 Using this last result, we shall now come back to the expression of $A^{-6} \Delta \beta$.
 First, we have
 \begin{eqnarray*}
 A^{- 2} &= &1 + 2 w \cos \alpha \cos t + w^2 \cos^2 t \\
  B & = & k + \left( k \cos \alpha + \frac{\sin \alpha \cos \alpha}{r} \right)
  w \cos t + \cos \alpha \left( k + \frac{2 \sin \alpha}{r} \right) w \cos t \\
	&& +
  \frac{\sin \alpha}{r} w^2 \cos^2 t\\
  & = & k + \left( 2 k \cos \alpha + \frac{3 \sin \alpha \cos \alpha}{r}
  \right) w \cos t + \frac{\sin \alpha}{r} w^2 \cos^2 t \\
  B' & = & k' + \left( 2 k' \cos \alpha - \frac{2 k^2}{\sin
  \alpha} - \frac{k}{r} (3 \cos^2 \alpha + 1) + \frac{\sin \alpha}{r^2} (3 -
  15 \cos^2 \alpha) \right) w \cos t\\
  &  & - \left( \frac{k \cos \alpha}{r} + \frac{6 \sin \alpha \cos
  \alpha}{r^2} \right) w^2 \cos^2 t
\end{eqnarray*}
After a
long but straightforward calculation, we get

\begin{eqnarray*}
   A^{-6} \Delta \beta & =   & -\left( k' + \frac{2 k \cos \alpha}{r} +
  \frac{3 kw^2 \cos \alpha}{r} + \frac{3 w^2 \sin \alpha \cos \alpha}{r^2} -
  \frac{\sin \alpha \cos \alpha}{r^2} \right)\\
   & &  - \left( 4 k' \cos \alpha + \frac{k^2}{\sin \alpha} + \frac{k (14
  \cos^2 \alpha - 2)}{r} - \frac{\sin \alpha}{r^2} (13 \cos^2 \alpha - 4) \right.\\
  && \hspace{1em} \left. +\frac{kw^2}{r} (6 \cos^2 \alpha + 3) + \frac{w^2 \sin \alpha (15 \cos^2
  \alpha)}{r^2} \right) w \cos t \\
  &  & - \left( k' (1 + 4 \cos^2 \alpha) + 5 k^2 \cot \alpha + \frac{k
  \cos \alpha}{r} (20 \cos^2 \alpha + 13) \right. \\
  && \quad \left. + \frac{\sin \alpha \cos
  \alpha}{r^2} (5 \cos^2 \alpha - 2) + \frac{6 k \cos \alpha w^2}{r} + \frac{6
  w^2 \sin \alpha \cos \alpha}{r^2} \right) w^2 \cos^2 t\\
  &  & - \left( 2 k' \cos \alpha + \frac{k^2}{\sin \alpha} (4 \cos^2
  \alpha) + \frac{k}{r} (22 \cos^2 \alpha + 1) \right. \\
  & & \quad \left.+ \frac{\sin \alpha}{r^2}
  \left( 7 \cos^2 \alpha + 4 + 3 w^2 \right) \right) w^3 \cos^3 t
\end{eqnarray*}
Thus the  vanishing of $\Delta \beta$ implies four linear relations in
the (dependent) variables $k' \cos \alpha$, $k^2 / \sin
\alpha$, $k / r$ and $\sin \alpha / r^2$, that we write in matrix form:
 $$ \left(\begin{array}{cccc}
     1 & 0 & X (2 + 3 w^2) & X (3 w^2 - 1)\\
     4 & 1 & 14 X - 2 + (6 X + 3) w^2 & 15 X w^2 + 4 - 13 X\\
     1 + 4 X & 5 X & X (20 X + 13) + 6 X w^2 & X(6 w^2+5 X - 2)\\
     2 & 4 X & 22 X + 1 & 3 w^2 + 7 X + 4
   \end{array}\right) \left( \begin{array}{c}
  k' \cos \alpha \\ k^2 / \sin \alpha \\ k / r \\ \sin \alpha / r^2
   \end{array}\right) =0$$
where we have denoted $X= \cos^2 \alpha$ for brevity. We reduce it to
\begin{equation} \label{system4x4}
\left(\begin{array}{cccc}
     1 & 0 & X (2 + 3 w^2) & X (3 w^2 - 1)\\
     0 & 1 & 6 X - 2 + (- 6 X + 3) w^2 & 3 X w^2 + 4 - 9 X\\
     0 & 0 & X (- 18 X + 21) + X (18 X - 12) w^2 & X (3 - 27 X)
     w^2 + X (54 X - 21)\\
     0 & 0 & - 24 X^2 + 26 X + 1 - X (24 X - 6) w^2 & w^2 (3 - 6 X -
     12 X^2) + 4 - 7 X + 36 X^2
   \end{array}\right) \end{equation}
leaving us with the following system
\begin{equation} \label{system2x2}
\left\{ \begin{array}{l}
     X( (- 18 X + 21) + (18 X - 12) w^2) k \hspace{3cm} \\
	\hspace{3cm} = X((27 X-3) w^2 + (21- 54 X))  \frac{\sin \alpha}{r}\\
     (24 X^2 - 26 X - 1 + X w^2 (24 X - 6)) k \hspace{3cm} \\
	\hspace{3cm} = -( w^2 (3 - 6 X - 12 X^2) + 4 - 7 X + 36 X^2) \frac{\sin \alpha}{r}
   \end{array} \right.
\end{equation}
The trivial solution $k = \sin \alpha = 0$ being excluded because $w$ is not defined when 
$\sin \alpha=0$, the determinant of the above system must vanish, hence the following algebraic equation
\begin{multline} 
E (X, Y) = \left( - 288 X^3 + 90 X^2 + 36 X - 12 \right) Y^2 \\
+ \left( - 504
   X^4 + 756 X^3 - 264 X^2 + 8 X + 4 \right) Y \\
+ 216 X^5 - 774 X^4 + 991 X^3 -
   489 X^2 + 21 X + 35 
\end{multline}
in terms of the variables $X = \cos^2 \alpha$ and $Y = c^2 r^{- 6}=w^2 \sin^2\alpha$, where we
have divided by $X$ for simplicity. We will now show a contradiction. Indeed
any of the two equations in (\ref{system2x2}) yields a dynamical system
\[ \left\{ \begin{array}{l}
     \alpha' = k - \frac{\sin \alpha}{r} = \frac{\sin \alpha}{r} g (\cos^2 \alpha, c^2 r^{- 6}) \\
     r' = \cos \alpha
   \end{array} \right. \]
and if we take for instance the first one
\[ g (X, Y) = \frac{(5 + 19 X + 12 X^2) (1 - X) + Y (3 - 36 X)}{(24 X^2 - 26 X
   - 1) (1 - X) + Y (24 X^2 - 6 X)} . \]
The solution $(\alpha (s), r (s))$ has to satisfy the analytic equation
$E (\cos^2 \alpha, c^2 r^{- 6}) = 0$. Taking the derivative leads to 
\[
2 \sin \alpha \cos \alpha \alpha' \partial E / \partial X + 6 c^2
r' r^{- 7} \partial E / \partial Y = 0, 
\]
which expands to
\[ F(X,Y) := (1 - X) g (X, Y) \frac{\partial E}{\partial X} + 3 Y \frac{\partial
   E}{\partial Y} = 0 \, .\]
For a solution to exist,
we need to find a common open set for the two algebraic curves $E = 0$
and $F = 0$. However $E$ is irreducible, so if $E$ and $F$ agree, they have to
do so on the complete component $E = 0$. We show on the contrary that this
component contains a point that does not satisfy $F = 0$. Let us solve $E (0,
Y) = 0$.
\[ 0 = E (0, Y) = - 12 Y^2 + 4 Y + 35 \]
yields $Y = \frac{1 \pm \sqrt{106}}{6}$. Evaluating $F$ at these two values
yields $1100 \pm 85 \sqrt{106}$ which are non zero, hence the contradiction. Finally we study
the case of $X = \cos^2\alpha=0$, where the system (\ref{system4x4}) amounts to
\[ \left\{ \begin{array}{l} \displaystyle
     \frac{k^2}{\sin \alpha} - \frac{2 k}{r} + \frac{3 kw^2}{r} + \frac{4 \sin \alpha}{r^2} = 0\\
\displaystyle     \frac{k}{r} + \frac{3 w^2 \sin \alpha}{r^2} + \frac{4 \sin \alpha}{r^2} 0
   \end{array} \right. \]
and since $k = \alpha' + \sin \alpha / r = \sin \alpha/r$
\[ \left\{ \begin{array}{l}
\displaystyle     \frac{\sin \alpha}{r^2} - \frac{2 \sin \alpha}{r^2} + \frac{3 \sin \alpha w^2}{r^2} 
+ \frac{4 \sin \alpha}{r^2} = 0\\
\displaystyle     \frac{\sin \alpha}{r^2} + \frac{3 w^2 \sin \alpha}{r^2} + \frac{4 \sin \alpha}{r^2} = 0
   \end{array} \right. \]
\[ 3 + 3 w^2 = 5 + 3 w^2 = 0 \]
so again we get a contradiction and the proof is complete.

\subsection{Hamiltonian stationary Type II surfaces}

\begin{proposition}
 A Hamiltonian stationary Lagrangian cyclic surface of $\C^2$ of type II
  is locally congruent to a Cartesian product of two round circles
$ \S^1(r_1) \times \S^1(r_2)$.
\end{proposition}

\noindent \textit{Proof.}
As we saw in Section 2, Remark~\ref{autrebase}, a type II immersion which is not a product of curves
can take the following form
$$\begin{array}{clcl} X: & I \times \R / 2 \pi \Z &  \to &  \C^2 \\
                    & (s,t) & \mapsto &
 c(\gamma_1(s) e^{i t}, \gamma_2(s) e^{i t}) +
  \int_{s_0}^s W(u)(\gamma'_1,\gamma'_2) du
      \end{array}$$
where the Legendrian curve $\gamma =
(\gamma_1, \gamma_2)$ is parametrized by arc length. We denote by $\beta_L =
\arg (\gamma_1 \gamma'_2 - \gamma_2 \gamma'_1)$ the Legendrian angle of the
curve $\gamma = (\gamma_1, \gamma_2)$, i.e. the Lagrangian angle of the plane
spanned by $\gamma, \gamma'$. Next we compute the derivatives of the immersion (we
set the constant $c$ to $1$ since that does not change the Lagrangian angle):
\[ Y_s = (W + e^{i t}) \left(\begin{array}{c}
     \gamma'_1\\
     \gamma'_2
   \end{array}\right) \;, \hspace{1em} Y_t = i e^{i t}  \left(\begin{array}{c}
     \gamma_1\\
     \gamma_2
   \end{array}\right), \]
from which we deduce that
\[ {\det}_{\C} (Y_s, Y_t) = i e^{i t} (W + e^{i t}) (\gamma'_1 \gamma_2 -
   \gamma'_2 \gamma_1) = - i (e^{i t} + W) e^{i (\beta_L + t)} \]
Thus
\[ \beta (s, t) = t + \beta_L + \arctan \left( \frac{\sin t + W_2}{\cos t +
   W_1} \right) - \frac{\pi}{2} \]
We now compute the induced metric:
\[ g_{ss} = 1 + |W|^2 + 2 \langle W, e^{i t} \rangle = : A^2 \quad g_{tt} = 1
   \quad g_{st} = 0, \]
and the first derivatives of $\beta :$
\[ \beta_s = \beta'_L + \frac{\langle W', i e^{i t} \rangle + W_1 W'_2 - W'_1
   W_2}{A^2} \;, \hspace{1em} \beta_t = \frac{2 + |W|^2 + 3 \langle W, e^{i t}
   \rangle}{A^2} \;, \]
\[  \]
where $\left\langle, \right\rangle$ stands for the scalar product in
$\C \simeq \R^2$.
\[ g = \left(\begin{array}{cc}
     A^2 & 0\\
     0 & 1
   \end{array}\right), \hspace{1em} g^{- 1} = \left(\begin{array}{cc}
     A^{- 2} & 0 \\
     0 & 1
   \end{array}\right) \; . \]
Thus,
\[ \Delta \beta = - \frac{1}{\sqrt{\det g}} \left( \frac{\partial}{\partial s}
   \left( \sqrt{\det g} g^{ss} \beta_s \right) + \frac{\partial}{\partial t}
   \left( \sqrt{\det g} g^{tt} \beta_t \right) \right) \]
\begin{eqnarray*}
  - A \Delta \beta & = & \frac{\partial}{\partial s} \left( \frac{1}{A}
  \left( \beta'_L + \frac{\langle W', i e^{i t} \rangle + W_1 W'_2 - W'_1
  W_2}{A^2} \right) \right) \\
	&& + \frac{\partial}{\partial t} \left( \frac{2 +
  |W|^2 + 3 \langle W, e^{i t} \rangle}{A} \right)
\\
  & = & \frac{\beta''_L A - \beta_L'  \frac{\partial A}{\partial s}}{A^2} +
  \frac{\left\langle W'', i e^{i t} \right\rangle + W_1 W''_2 - W''_1 W_2
  }{A^3}\\
  &  & - 3 \frac{\langle W', i e^{i t} \rangle + W_1 W'_2 - W'_1 W_2}{A^4}
  \frac{\partial A}{\partial s}\\
  &  & + 3 \frac{\left\langle W, i e^{i t} \right\rangle}{A} - \frac{2 + |W|^2
  + 3 \langle W, e^{i t} \rangle}{A^2}  \frac{\partial A}{\partial t}\\
  & = & \frac{\beta''_L}{A} - \frac{\beta_L'  \left\langle W', W + e^{i t}
  \right\rangle}{A^3} + \frac{\left\langle W'', i e^{i t} \right\rangle + W_1
  W''_2 - W''_1 W_2 }{A^3}\\
  &  & - 3 \frac{\left( \langle W', i e^{i t} \rangle + W_1 W'_2 - W'_1 W_2
  \right)  \left\langle W', W + e^{i t} \right\rangle}{A^5}\\
  &  & + 3 \frac{\left\langle W, i e^{i t} \right\rangle}{A} - \frac{2 + |W|^2
  + 3 \langle W, e^{i t} \rangle}{A^3}  \left\langle W, i e^{i t} \right\rangle
\end{eqnarray*}
\begin{eqnarray*}
  - A^6 \Delta \beta & = & A^4 \beta''_L - A^2 \beta_L'  \left\langle W', W +
  e^{i t} \right\rangle + A^2 ( \left\langle W'', i e^{i t} \right\rangle + W_1
  W''_2 - W''_1 W_2)\\
  &  & - 3 \left( \langle W', i e^{i t} \rangle + W_1 W'_2 - W'_1 W_2 \right)
  \left\langle W', W + e^{i t} \right\rangle\\
  &  & + A^2  \left\langle W, i e^{i t} \right\rangle  \left( 1 + 2| W|^2 + 3
  \left\langle W, e^{i t} \right\rangle  \right) .
\end{eqnarray*}
For fixed $s$, the expression $- A^6 \Delta \beta$ is polynomial of degree 3
in the variable $e^{i t}$, and the only term of order 3 comes from the last
line in the expression above: $3 A^2  \left\langle W, e^{i t} \right\rangle
\left\langle W, i e^{i t} \right\rangle$. Thus if $\beta$ is harmonic $W$ must
vanish, hence we are in the centered case (Hopf type tori).

This case is straightforward: on the one hand the metric is flat and on the
other hand $\beta = 2 t + \beta_L - \frac{\pi}{2}$, thus we must have
$\beta''_L = 0$. Following the terminology of [CLU], the curve $\gamma$ is
\textit{contact stationary}
 and take the following form :
\[ \gamma (s) =  (
   c e^{is/\sqrt{2}c }, \sqrt{1-c^2} e^{-is/\sqrt{2(1-c^2)}}), \]
where $c \in (0,1)$. We observe that the closedness of the Legendrian curve is not necessary in order
to get a compact surface. The corresponding
Hopf torus is a product of circles $\S^1(c) \times \S^1(\sqrt{1-c^2})$.
By making a convenient homothety we get any product of circles $\S^1(r_1) \times \S^1(r_2)$ thus the proof
is complete.

\section{Ruled Lagrangian surfaces}
Let $\Sigma$ be a smooth, ruled, Lagrangian surface.
If the rulings are parallel, it is straightforward that $\Sigma$
is a Cartesian product ${\cal L} \times {\Gamma}$ of some straight line ${\cal L}$,
and some planar curve $\Gamma$,
such that ${\cal L} \subset P_1$ and ${\Gamma} \subset P_2$, where $P_1$ and $P_2$ are two
orthogonal, complex planes. So from now on, our discussion being local, we shall assume
that the rulings are not parallel.

Locally, $\Sigma$ may be parametrized by the following immersion:
 $$\begin{array}{clcl} X: & I \times \R  &  \to &  \C^2 \\
                    & (s,t) & \mapsto &  \gamma(s) t + V(s),
      \end{array}$$
where $\gamma(s)$ is a
unit speed curve of $\S^3$ and $V(s) \in \C^2$,
 and we have
$$X_s= \gamma' t + V', \hspace{3em} X_t = \gamma.$$

We now claim that, without loss of generality, we may assume that
$X_s$ and $X_t$ are orthogonal: to see this, we reparametrize the surface by
$\hat{X}(s,t) = X(s,t) + \gamma(s) T(s)$, where $T(s)$ is some real function. We observe
that the images of $X$ and $\hat{X}$ are the same. Then we compute
$$ \hat{X}_s=X_s + \gamma' T + \gamma T',
  \hspace{3em}  \hat{X}_t=X_t=\gamma .$$
Thus by choosing $T$ such that $T' = - \<X_s, X_t\>$, we may reparametrize our
surface such that $\hat{X}_s$ and $\hat{X}_t$ are orthogonal.

The Lagrangian assumption amounts to
 $$ \omega(X_s,X_t) = \<X_s, J X_t \>
               = t \< \gamma', J \gamma \> + \< V', J\gamma \> =0.$$
 The vanishing of $\< \gamma', J \gamma \>$ means that the curve $\gamma$ is Legendrian. It follows
 that $(\gamma , J \gamma , \gamma', J \gamma')$ is an orthonormal basis of $\R^4 \simeq \C^2$.
Thus the conditions $\<V', J\gamma \> =0$ and $\<X_s,X_t\>=\<V', \gamma\>= 0$ imply that
there exists a planar curve
$\alpha(s)=x(s)+iy(s)$ such that
 $V'= x \gamma' + y J \gamma'= \alpha  \gamma'$. So we have shown the first part of the following:

\begin{theorem}
A smooth, ruled, Lagrangian surface of $\C^2$ is either the Cartesian
 product $\mathcal{L} \times \Gamma$ of a  straight line $\mathcal{L}$ with a planar curve $\Gamma$  or  may be locally parametrized by an immersion of the
form
$$\begin{array}{clcl} X: & I \times \R  &  \to &  \C^2 \\
                    & (s,t) & \mapsto &  \gamma(s) t +
           \int_{s_0}^s \alpha(u) (\gamma_1'(u),\gamma_2'(u))du,
      \end{array}$$
where $\gamma=(\gamma_1,\gamma_2)$ is some Legendrian curve of $\S^3$ and $\alpha=(x,y)$ is some
planar curve. Moreover, the only self-similar ruled Lagrangian surfaces of $\C^2$ are
Cartesian products $\mathcal{L} \times \Gamma$, where $\mathcal{L}$ is a straight line and
${\Gamma}$ is a planar self-similar  curve (cf \cite{AL}).
\end{theorem}

\noindent \textit{Proof.}
We have already characterized ruled Lagrangian surfaces, so it remains to study the
self-similar equation.

We start computing the first derivatives of the immersion:
$$ X_t(s,t)=\gamma(s), \hspace{3em}
 X_s(s,t)=
(1 + \alpha(s))(\gamma'_1(s),\gamma'_2(s))$$
 from which we deduce the expression of the induced metric:
$$ E=(t+x)^2 + y^2, \hspace{3em} G= 1, \hspace{3em} F=0,$$
and a basis of the normal bundle:
$$ N_t= i \gamma(s), \hspace{3em}
 N_s= i(1 + \alpha(s))(\gamma'_1(s),\gamma'_2(s)).$$
We now compute some second derivatives:
$$ X_{st}= \gamma'(s), \hspace{3em} X_{tt}= 0.$$
We deduce that
$$ 2\< H, N_t \> = \frac{\<X_{ss},N_t\>}{E} + \frac{\<X_{tt},N_t\>}{G}
  =\frac{\<X_{st},N_s\>}{E}= \frac{ 1 + x}{(t+x)^2 + y^2}.$$

On another hand, we have:
$$ \<X, N_t \>=\<i \gamma(s), \int_{s_0}^s \alpha(u) (\gamma_1'(u),\gamma_2'(u))du \>,$$
which clearly does not depend on $t$.
 So the equation
$$ \<H, N_t\> + \lambda \<X, N_t \>=0$$
 can never hold for such an immersion, so the only self-similar ruled Lagrangian surfaces
 are products ${\cal L} \times {\Gamma}$. It is then straightforward to see that the curve ${\Gamma}$
must be a solution of the self-similar equation.

\begin{remark} We recover the Lagrangian helicoid of \cite{B} by taking
$\gamma(s)=(k+il)(\cos s, \sin s)$, where $k$ and $l$ are two constants such
that $k^2+l^2=1$. In the notation of \cite{B}, we have $x(s)=G/2$ and $y(s)=-A$.

\end{remark}

\bigskip

\noindent
Henri Anciaux \\
Departement of Mathematics and Computing \\
Institute of Technology, Tralee \\
Co. Kerry, Ireland \\
henri.anciaux@staff.ittralee.ie \\

\medskip

\noindent
Pascal Romon \\
Universit\'{e} de Paris-Est Marne-la-Vall\'{e}e \\
5, bd Descartes, Champs-sur-Marne \\
77454 Marne-la-Vall\'{e}e cedex 2, France \\
pascal.romon@univ-mlv.fr

\end{document}